\documentclass{amsart}
\usepackage{amsmath,amscd,amssymb,amsthm,epsfig,psfrag,latexsym,subfigure}

\newtheorem{theorem}{Theorem}[section]
\newtheorem{lemma}[theorem]{Lemma}
\newtheorem{corollary}[theorem]{Corollary}

\newtheorem{definition}[theorem]{Definition}

\theoremstyle{remark}
\newtheorem*{acknowledgments}{Acknowledgments}

\def\vol{\mbox{\rm{Vol}}}

\def\b{\beta}

\def\g{\gamma}

\def\S{\Sigma}

\def\int{\mbox{\rm{int}}}

\def\N{\mathcal{N}}

\begin{document}
\title{Small 3-manifolds of large genus}
\author{Ian Agol}
\address{Department of Mathematics, University of Illinois at
Chicago, 322 SEO m/c 249, 851 S. Morgan Street, Chicago, IL
60607-7045}

\email{agol@math.uic.edu}
\thanks{Partially supported by ARC grant 420998}


\begin{abstract}
We prove the existence of pure braids with arbitrarily many
strands which are small, {\it i.e.} they contain no closed
incompressible surface in the complement which is not boundary
parallel. This implies the existence of irreducible non-Haken
3-manifolds of arbitrarily high Heegaard genus.

\end{abstract}
\date{\today}
\maketitle

\section{Introduction}
Haken introduced the notion of irreducible 3-manifolds containing
an incompressible surface, called sufficiently large, or Haken
\cite{Hak68}. These were thrust into prominence by theorems of
Waldhausen, who (among other things) showed that the homeomorphism
and word problems are solvable for these manifolds \cite{Wal:78},
and by Thurston who showed that they satisfy the geometrization
conjecture \cite{Th:82}. It is therefore interesting to understand
how prevalent non-Haken 3-manifolds are. Jaco and Shalen showed
that a Seifert fibred space $M$ is non-Haken if and only if it is
atoroidal with $H_1(M)$ finite \cite{JS79}. Thurston showed that
all but finitely many Dehn fillings on the figure eight knot
complement are non-Haken \cite{Th}. Hatcher and Thurston extended
this result to all 2-bridge knot exteriors \cite{HT85}. These
examples all have Heegaard genus 2. By a result of Hatcher, if a
link is small, {\it i.e.} is irreducible and contains no closed
incompressible surface other than boundary tori, then Dehn filling
on any boundary component yields a small manifold for all but
finitely many fillings \cite{Ha82}. Floyd and Hatcher \cite{FH82}
and independently Culler, Jaco, and Rubinstein \cite{CJR82} showed
that punctured torus bundles and 4-punctured sphere bundles with
irreducible monodromy are small. Oertel proved that Montesinos
knots covered by small Seifert fibred spaces are small \cite{O84}.
It is believed that the Seifert-Weber dodecahedral space is small.
Several other examples of small manifolds are known, {\it e.g.}
the Borromean rings (Lozano), some chain links with $\leq 5$
components (Oertel), and some other examples of Dunfield
\cite{D99}, Hass and Menasco \cite{HM93}, and Lopez
\cite{Lo92,Lo93}.

A result of Moriah and Rubinstein shows that the Heegaard genus of
infinitely many Dehn fillings on a link complement is the same as
that of the link complement \cite{MR97}. Since most punctured
torus bundles have Heegaard genus 3 \cite{CS99}, there are
infinitely many closed small 3-manifolds of genus 3. This appears
to be the largest known genus of a small 3-manifold in the
literature. Alan Reid asked whether there are small links of
arbitrarily many components, observing that this would imply the
existence of irreducible non-Haken 3-manifolds of arbitrarily
large Heegaard genus, as we show in theorem \ref{small}. We prove
in theorem \ref{smallbraids} the existence of small links of
arbitrarily many components, answering Reid's question.

\begin{acknowledgments}
We thank Alan Reid, Hyam Rubinstein, Saul Schleimer, and Bill
Thurston for contributing ideas and suggestions which led to the
results in this paper. The ideas in the paper were partly inspired
by a talk of Dan Margalit at UIC on the pants complex
\cite{Mar02}, and by the work of Jeff Brock \cite{Bro01}.
\end{acknowledgments}
\section{Pants decompositions}
\begin{figure}[htb]
    \begin{center}
    \subfigure{\epsfig{figure=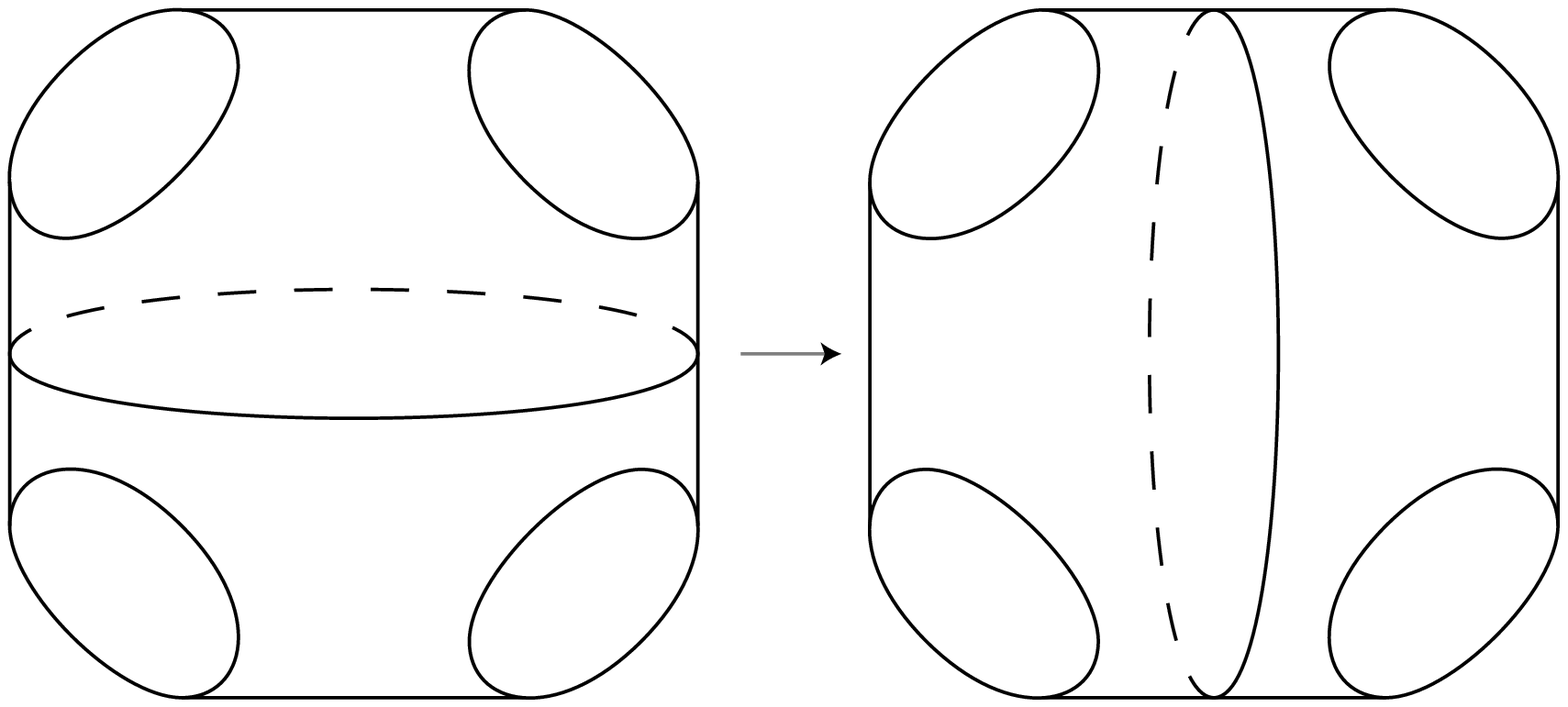,width=.45\textwidth}}\qquad
    \subfigure{\epsfig{figure=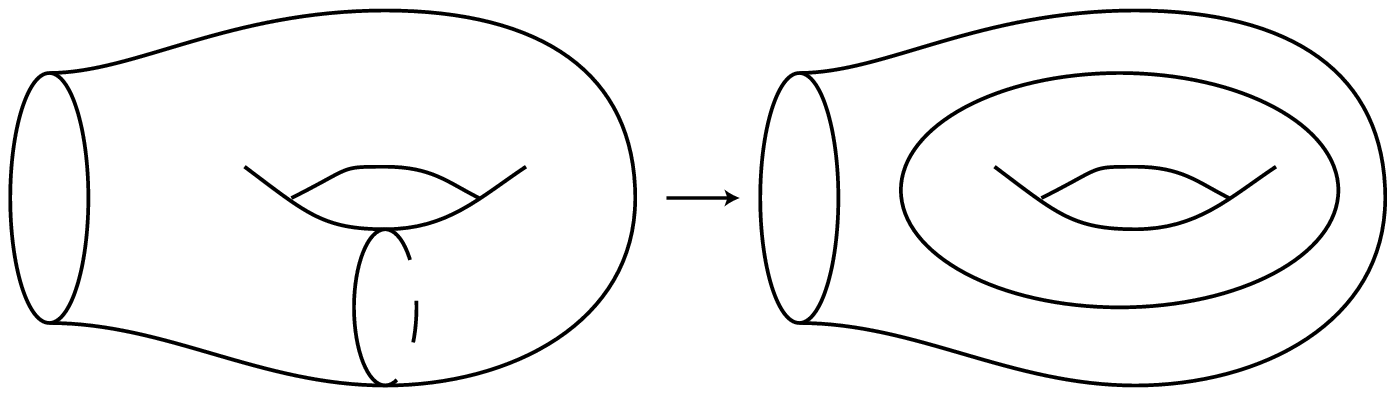,width=.45\textwidth}}
    \caption{An $A$-move and $S$-move \label{Amove}}
    \end{center}
\end{figure}

The pants complex of a surface was introduced by Hatcher, Lochak,
and Schneps \cite{HLS00}, and we will follow their conventions and
terminology . We will just be interested in the pants graph. Let
$\Sigma$ be a connected compact orientable surface with
$\chi(\Sigma)<0$ . We say $\Sigma$ has type $(g,n)$ if it has
genus $g$ and $n$ boundary components. By a {\bf pants
decomposition} of $\Sigma$, we mean a finite collection $P$ of
disjoint smoothly embedded circles cutting $\Sigma$ into pieces
which are surfaces of type $(0,3)$. The number of curves in a
pants decomposition is $3g-3+n$, and the number of pants is
$-\chi(\Sigma)$. Let $P$ be a pants decomposition, and suppose
that one of the circles $\beta$ of $P$ is such that deleting
$\beta$ from $P$ produces a complementary component of type
$(1,1)$. This is equivalent to saying there is a circle $\g$ in
$\Sigma$ which intersects $\b$ in one point transversely and is
disjoint from all the other circles in $P$. In this case,
replacing $\b$ by $\g$ in $P$ produces a new pants decomposition
$P'$ (fig. \ref{Amove}).
 We call this replacement a {\bf simple move}, or
$S$-move. In a similar fashion, if $P$ contains a circle $\b$ such
that deleting $\b$ from $P$ produces a complementary component of
type $(0,4)$, then we obtain a new pants decomposition $P'$ by
replacing $\b$ with a curve $\g$ intersecting $\b$ transversely in
two points and disjoint from the other curves of $P$ (fig.
\ref{Amove}). The transformation $P\rightarrow P'$ in this case is
called an {\bf associativity move}, or $A$-move.

\begin{definition}
The {\bf pants decomposition graph} $\mathcal{P}(\Sigma)^{(1)}$ is
the graph having vertices the isotopy classes of pants
decompositions of $\S$, with an edge joining two vertices whenever
the corresponding pants decompositions differ by a single $A$- or
$S$-move. This is the one skeleton of the pants decomposition
complex, which was proven to be connected in \cite{HLS00}.
\end{definition}

{\bf Main Construction:} Suppose we have a homeomorphism
$\psi:\S\longrightarrow \S$, with mapping torus $T_{\psi}=(S\times
I)/\{(x,0)\sim (\psi(x),1)\}$ . Given a path
$C=P_0-P_1-\cdots-P_{m-1}-P_m \subset \mathcal {P}(\Sigma)^{(1)}$,
$P_m=\psi(P_0)$,  we get a sequence of circles $\b_1,...,\b_m$,
where $\b_{i+1}$ is the circle in $P_{i+1}$ replacing a circle in
$P_i$ ($i$ is taken $\pmod{m}$). We assume that there is no circle
$\b_j$ which is contained in all of the partitions $P_i$, $1\leq
i\leq m$. For each circle $\b_i$, drill out a curve
$B_i=\b_i\times \{\frac{i}{m}\}$ in $T_\psi$. For each
complementary region of $P_i$ with boundary curves
$\b_i,\b_j,\b_k$ (where $i,j,k$ might not be distinct), there is
an embedded pants in $T_{\psi}$ with boundary on the curves
$B_i\cup B_j\cup B_k$ and interior disjoint from $\cup_i B_i$, and
moreover, all of these pants may be chosen to have disjoint
interiors (we pull apart the pants lying in $\S$ like an
accordion). For an $A$- or $S$-move $P_i-P_{i+1}$, there is a
complementary region bounded by two pants for an $S$-move, and
four pants for an $A$-move, which we call an $A$- or $S$-region,
respectively. Consider the link complement $M_C=T_{\psi}\backslash
\N(\cup_i B_i)$. Then $T_{\psi}$ is obtained from $M_C$ by Dehn
filling on the boundary components corresponding to the link
$\cup_i B_i$. $M_C$ is decomposed along pants into $A$- or
$S$-regions.

\begin{definition}
A surface  $(S,\partial S) \subset (M_C,\partial M_C)$ (where if
$C=\emptyset$, then $M_C=T_\psi$) is {\bf pairwise incompressible}
if

\begin{enumerate}
\item
each component of $\partial S$ is either parallel to a component
of  $\partial \S$ in $\partial T_\psi$ or to a longitude of a
component of $\partial \N(B_i)$,
\item
if there is an annulus $A\subset M_C$ such that  $\int A\cap
S=\emptyset$, one boundary component of $A$ is in $\int S$ and the
other is in $\partial M_C$ parallel to a component of $\partial
\S$ or a longitude of $\partial \N(B_i)$, then $A$ is isotopic
into $S$.
\end{enumerate}
\end{definition}
\begin{lemma} \label{tube}
$\int M_C$ has a complete hyperbolic metric of finite volume. Any
pairwise incompressible surface in $M_C$ is a disjoint union of
pants. Any pairwise incompressible surface in $T_\psi$ isotopic
into $M_C$ is obtained by tubing pants in $M_C$.
\end{lemma}
\begin{proof}
The $A$- and $S$-regions may be regarded as pared manifolds
\cite{Mo}, with pared locus consisting of annuli in the boundary
which are regular neighborhoods of the curves from the pants
decomposition. Each $A$-region may be decomposed into two ideal
octahedra (fig. \ref{Aregion}), and each $S$-region is decomposed
into one ideal octahedron (fig. \ref{Sregion}).

\begin{figure}[htb]
    \begin{center}
    \epsfxsize= \textwidth
    \epsfbox{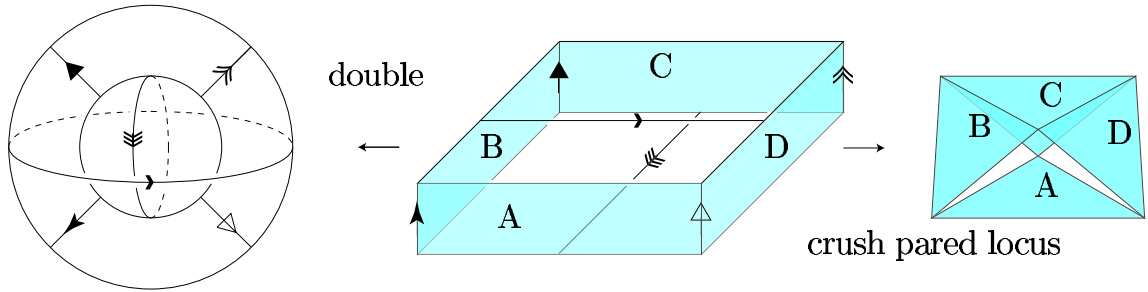}
    \caption{An $A$-region from 2 octahedra \label{Aregion}}
    \end{center}
\end{figure}

\begin{figure}[htb]
    \begin{center}
    \epsfxsize= \textwidth
    \epsfbox{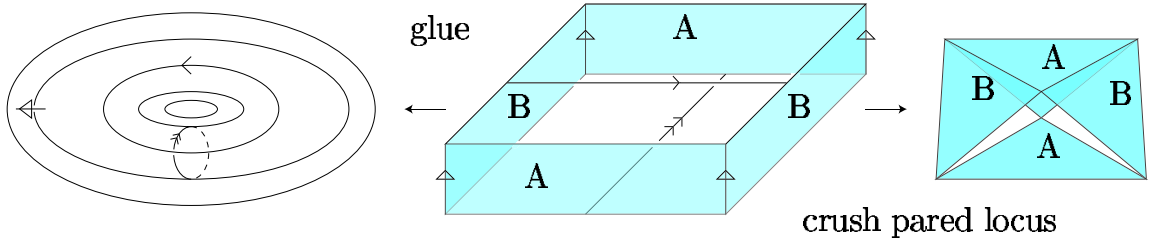}
    \caption{An $S$-region from one octahedron \label{Sregion}}
    \end{center}
\end{figure}

Thus, an $A$-region is obtained by doubling a checkerboard colored
ideal octahedron along the dark faces. An $S$-region is obtained
by folding pairs of dark faces of a checkerboard colored ideal
octahedron together along common ideal vertices. This gives each
$A$- and $S$-region a hyperbolic metric with totally geodesic
boundary, and rank one cusps along the pared locus. We may glue
these pieces together along geodesic 3-punctured spheres to get a
hyperbolic structure on $\int M_C$.

Suppose that we have an incompressible surface $(S,\partial
S)\subset  (M_C,\partial T_\psi)$, such that $S$ is pairwise
incompressible in $T_\psi$. We may assume that it intersects each
pants in essential simple closed curves. Since closed curves in
pants are boundary parallel, we may do surgery along annuli, to
get a surface $S'$ which is disjoint from the pants (fig.
\ref{surger}). Since $M_C$ is acylindrical, we may assume that
none of the surgeries produce annuli.
\begin{figure}[htb]
    \begin{center}
    \epsfxsize= 2 in
    \epsfbox{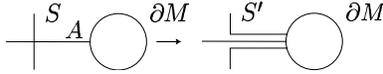}
    \caption{Surgering $S$ along an annulus $A$ to obtain $S'$ \label{surger}}
    \end{center}
\end{figure}
Thus, each component of the surface $S'$ lies in an $A$- or
$S$-region. Since an ideal octahedron is the same as a truncated
tetrahedron with its edges drilled out (where dark faces of the
octahedron correspond to faces of the tetrahedron), we may assume
that the surface intersects each octahedron of the decomposition
of $A$- and $S$-regions normally. A surface is normal if and only
if it is incompressible in the complement of the 1-skeleton
\cite{Tho}, so since there is no one skeleton, the surface must be
normal. In the case of an $A$-region, the surface must be a double
of a quadrilateral or a triangle. The doubles of triangles give
the pants surfaces, and the doubles of the quadrilaterals have
annular compressions to the pared locus, that is they are obtained
by tubing together pairs of pants. A similar property holds for
the $S$-regions. Thus, we may obtain every pairwise incompressible
surface by tubing together pants.
\end{proof}

\begin{corollary}
$\vol(T_\psi)\leq V_{oct}(2A+S)$, where $\vol$ denotes the volume
of the complete hyperbolic metric, $V_{oct}$ denotes the volume of
a regular ideal octahedron, $A$ is the number of $A$-moves in the
path $C$, and $S$ denotes the number of $S$-moves.
\end{corollary}
\begin{proof}
Each $A$-region of $M_C$ contributes $2 V_{oct}$, and each
$S$-region contributes $V_{oct}$. Then by a theorem of Thurston
\cite{Th}, $\vol(T_\psi)< \vol(M_C)=V_{oct}(2A+S)$.
\end{proof}
{\bf Remark:} This gives a tight upper bound for a theorem of
Brock \cite{Bro01} relating volumes of mapping tori to pants
distance of the monodromy.

\section{Small links}
\begin{theorem} \label{smallbraids}
There exist pure braids in $S^2\times S^1$ with arbitrarily many
components and no closed incompressible surface in the complement.
\end{theorem}

\begin{theorem} \label{small}
There are small 3-manifolds of arbitrarily large Heegaard genus.
\end{theorem}
\begin{proof}\ref{small}
The following argument was observed by Alan Reid. By theorem
\ref{smallbraids}, there is a small link with $n$ components, for
$n$ arbitarily large. A link complement with $n$ boundary
components has Heegaard genus $g\geq n/2$, since there must be
$\geq n/2$ boundary components lying to one side of a Heegaard
surface, and therefore the genus of the Heegaard surface is $\geq
n/2$. A result of Hatcher \cite{Ha82} shows that there are
finitely many Lagrangian subspaces of the space of measured
laminations on the boundary which consist of laminations which are
the boundary of incompressible measured laminations. Moriah and
Rubinstein \cite{MR97} prove that for each boundary component,
there are finitely many points and lines in Dehn filling space of
the component, such that if a Dehn filling avoids these, then the
Heegaard genus of the resulting manifold is the same as the
Heegaard genus of the link complement. Combining these two
theorems, we see that there are infinitely many Dehn fillings on
the link which are small with genus $\geq n/2$, by doing Dehn
fillings which avoid the slopes given by Moriah and Rubinstein's
theorem and avoiding the boundary slopes of surfaces given by
Hatcher's theorem.
\end{proof}

\begin{proof} \ref{smallbraids}
Let $\Sigma_n$ be the $n$-punctured sphere, {\it i.e.} a surface
of type $(0,n)$. We construct a closed path $C$ in the pants graph
$\mathcal{P}(\S_n)^{(1)}$ and apply the Main Construction to
obtain a sequence of loops in $\Sigma_n\times S^1$. Then we do
very high Dehn twists about these loops to force any
incompressible surface in the resulting braid to be isotopic into
the complement of these loops (using Hatcher's theorem). By lemma
\ref{tube}, any pairwise incompressible surface is obtained by
tubing together pants. Then we analyze the ways of tubing pants
together, and show that the resulting surfaces are always
compressible in the braid complement obtained by any Dehn twists
about these loops. Bill Thurston suggested considering manifolds
which fiber over $S^1$, Saul Schleimer suggested drilling out
horizontal curves in a braid to try to find small braids by large
Dehn twists, and Hyam Rubinstein suggested using links whose
complements decompose along pants.

The first observation is that if there is a closed incompressible
surface in $T_\psi$, then there is a pairwise incompressible
surface, obtained by doing all possible compressions along
compressing annuli.

 We choose a closed path in the pants graph
$\mathcal{P}(\S_n)^{(1)}$, such that each pants is a twice
punctured disk or a punctured annulus, when thought of as lying in
$\S_n\times S^1/\subset S^2\times S^1$ . If this is true, then
after tubing along the horizontal loops, the surface becomes a
punctured sphere, torus, or klein bottle. If there is a
2-punctured disk, then the surface must be a punctured sphere. The
only essential punctured sphere in a braid complement is the fiber
surface. To see this, we take a cyclic cover of $S^2\times S^1$,
and lift the punctured $S^2$ to this cover, so that it lies
between two fiber surfaces. But the only incompressible surfaces
in the complement of the fiber are copies of the fiber. If copies
of the fiber are tubed together, then the resulting surface is
compressible. Thus, we may assume that the only pants that occur
are punctured annuli.

Now we construct an explicit closed path in
$\mathcal{P}(\S_n)^{(1)}$. This path seems to be the simplest to
construct with the property that each pants in the path is a
punctured annulus or 2-punctured disk. Consider the punctures of
$\S_n$ as lying in a great circle $\gamma$ on $S^2$. Cyclically
label the gaps between punctures along $\g$ $1,\ldots,n$. For
every pair of gaps $(i,j)$, there is a loop $\b_{i,j}$ in $S^2$
which meets $\gamma$ precisely twice at the gaps $i$ and $j$. We
form a path $C \subset \mathcal{P}(\S_n)^{(1)}$. The initial pants
decomposition is given by $P_0=\{\b_{1,3},\b_{1,4},...,
\b_{1,n-1}\}$. We may take another pants decomposition
$P_{n-3}=\{\b_{2,4},\b_{2,5}, ..., \b_{2,n}\}$, which is obtained
from $P_0$ by shifting each index by 1. There is a path of pants
decompositions $C_0=P_0-P_1-\ldots-P_{n-3}$, where $P_k$ is
obtained by replacing the first $k$ loops of $P_0$ with the first
$k$ loops of $P_{n-3}$, $k=1,\ldots,n-3$ (fig. \ref{pantspath7}).
More generally, let
$P_{j(n-3)}=\{\b_{j+1,j+3},\b_{j+1,j+4},\ldots,\b_{j+1,j+n-1}\}$,
where indices are taken $\pmod{n}$. For $1\leq k\leq n-3$, let
$P_{j(n-3)+k}$ be obtained by replacing the first $k$ loops of
$P_{j(n-3)}$ with the first $k$ loops of $P_{(j+1)(n-3)}$. Then we
have a closed path $C=P_0-P_1-\cdots-P_{n(n-3)}$ (fig.
\ref{pantspath5}).

\begin{figure}[htb]
    \begin{center}
    \epsfxsize=\textwidth
    \epsfbox{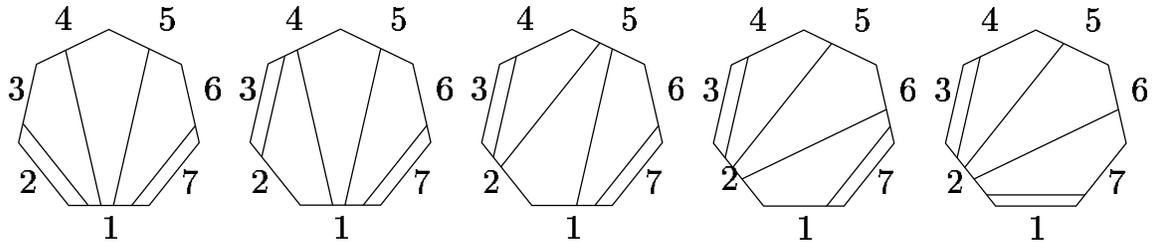}
    \caption{The path $C_0$ for $\S_7$. Double each polygon to obtain loops on $\S_7$,
    with punctures at vertices. \label{pantspath7}}
    \end{center}
\end{figure}

\begin{figure}[htb]
    \begin{center}
    \epsfxsize= \textwidth
    \epsfbox{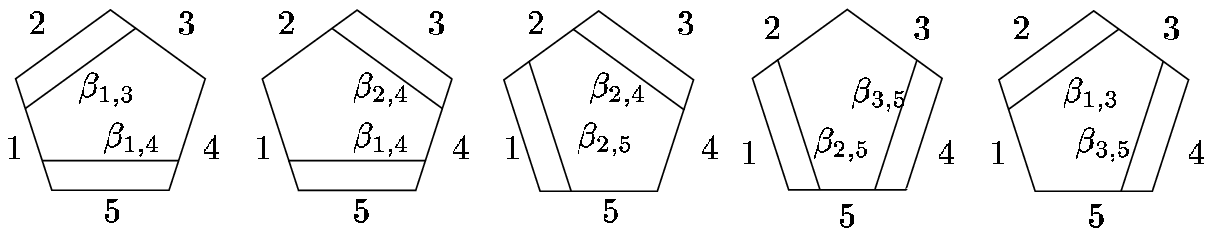}
    \caption{Half of the path $C$ for $\S_5$ \label{pantspath5}}
    \end{center}
\end{figure}

Any curve $\b_{i,j}$ such that $i-j\equiv \pm 2\pmod{n}$ bounds a
twice punctured disk in $\Sigma$. There are once punctured annuli
between the circles $\b_{i,j}$ and $\b_{i, j+1}$ or between
$\b_{i,j}$ and $\b_{i+1,j}$. Create loops $B_{i,j}$ in $\S_n\times
S^1$, each corresponding to the circle $\b_{i,j}$ in $\S_n$, as in
the Main Construction. Figures \ref{braid5}, \ref{braid6},
\ref{braid7} show a picture of part of $M_C$ for 5, 6, and 7
strand braids. The pictures are of a neighborhood of $\g\times I$,
cut along a gap, and flattened out. Loops $B_{i,j}$ and
$B_{i-1,j+1}$ have been isotoped to the same level.

\begin{figure}[htb]
    \begin{center}
    \epsfxsize= 3 in
    \epsfbox{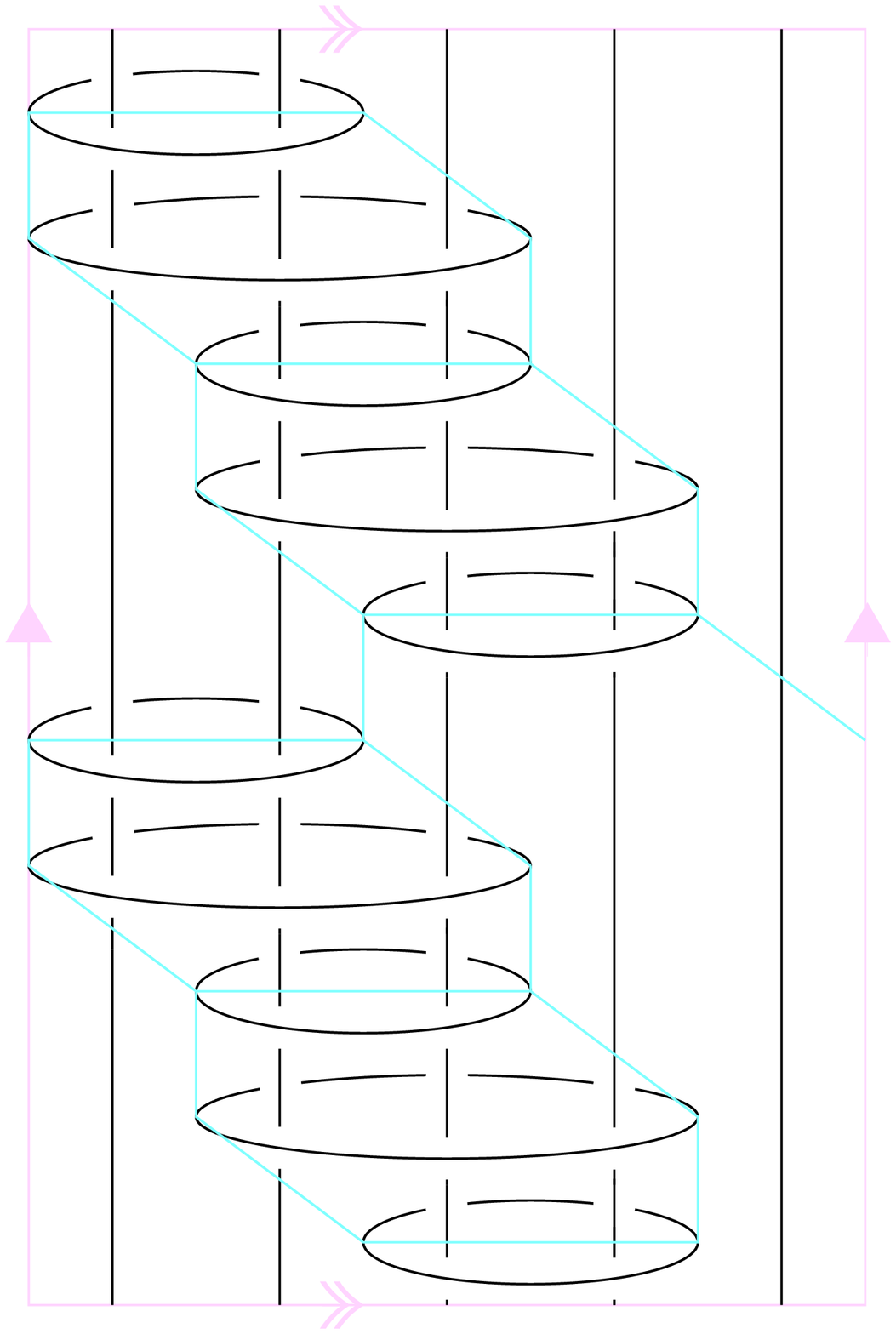}
    \caption{$M_C$ for $\S_5$ \label{braid5}}
    \end{center}
\end{figure}
\begin{figure}[htb]
    \begin{center}
    \epsfxsize= 3 in
    \epsfbox{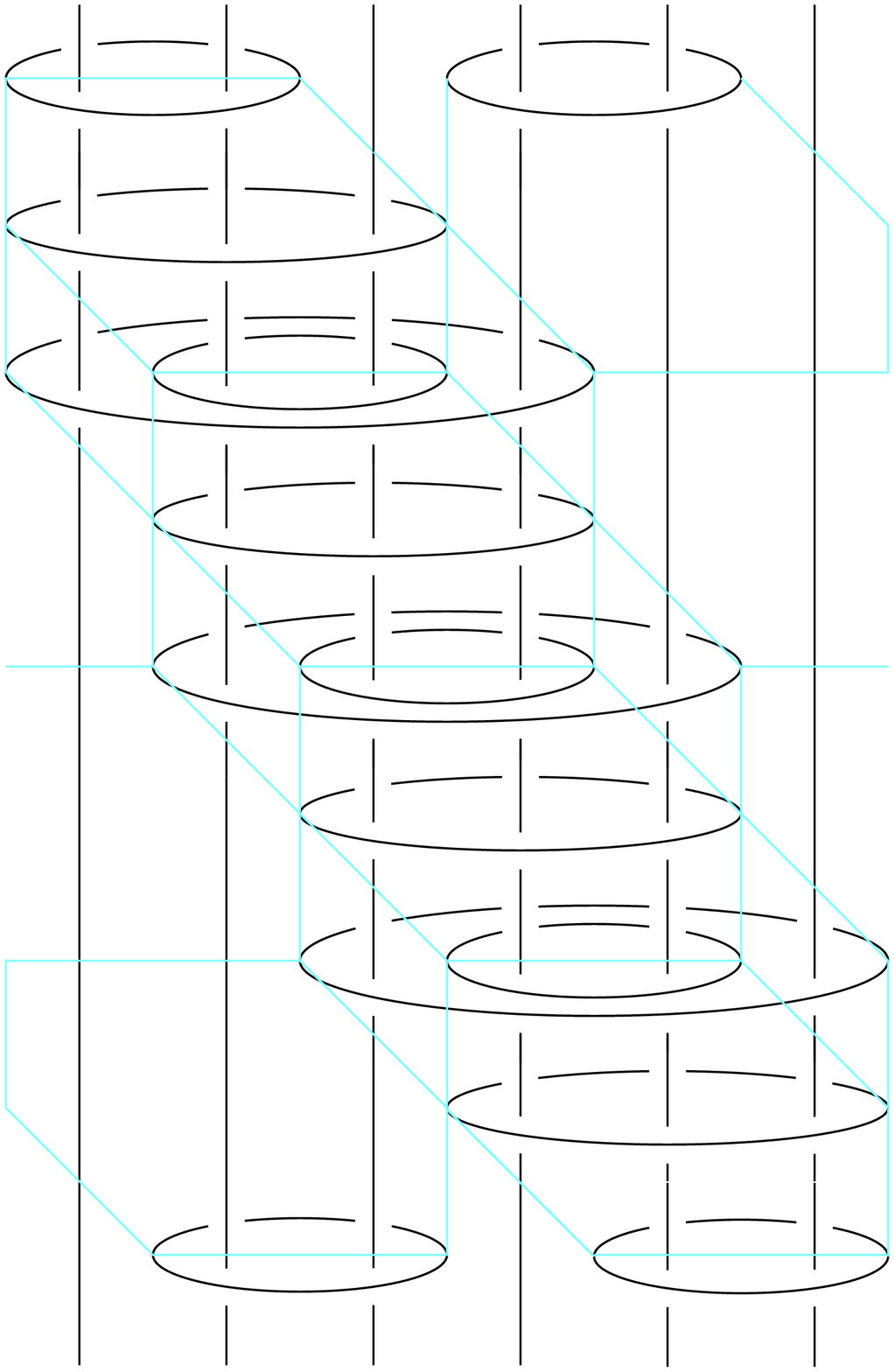}
    \caption{Part of $M_C$ for $\S_6$ \label{braid6}}
    \end{center}
\end{figure}
\begin{figure}[htb]
    \begin{center}
    \epsfxsize= 3 in
    \epsfbox{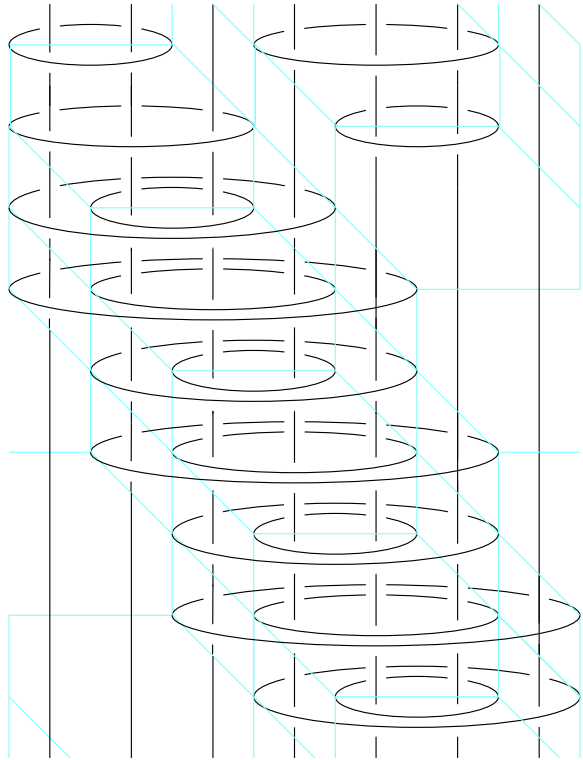}
    \caption{Part of $M_C$ for $\S_7$ \label{braid7}}
    \end{center}
\end{figure}

A key observation is that if one does Dehn twists about the
horizontal loops, {\it i.e.} by $1/n$ Dehn filling on these
curves, then an incompressible surface which is obtained by tubing
together pants in $M$ may be isotoped across the loops without
affecting its isotopy class in the Dehn filling. This follows
since each pants is isotopic into a fiber surface, so the framing
induced by the pants is the same as that of the fiber. For this
reason, we may draw surfaces in $M_C$ going through loops
$B_{i,j}$ without drawing which side of the loop the surface lies
on. Also, this makes it clear that two parallel pants may not be
tubed together to get an incompressible surface.

We have a 2-complex in the trivial braid which consists of the
loops and punctured disks and annuli (this is shown for $n=5,6,7$
in figures \ref{braid5},\ref{braid6}, \ref{braid7}) . The
intersection with $\g\times I$ is shown in figure \ref{complex}(a)
for $\S_7$. Take the subcomplex consisting only of punctured
annuli (fig. \ref{complex}(b)). Then one can see that any surface
carried by this complex is a punctured torus, isotopic to the
punctured torus $T$ which consists of the annuli bounded by the
sequence
${B_{1,3},B_{1,4},B_{2,4},B_{2,5},B_{3,5},...,B_{n,3},B_{1,3}}$
(fig. \ref{complex}(d)). If $n=5$, then this is the only
possibility for a surface. To see that every surface is isotopic
to $T$ if $n>5$, note that the punctured annulus spanning
$B_{i,j},B_{i,j+1},B_{i+1,j+1}$ is isotopic to the punctured
annulus spanning $B_{i,j},B_{i+1,j},B_{i+1,j+1}$ (this follows
from the symmetry of the $A$-regions: tubing two once punctured
annuli gives a 4-punctured sphere, which is isotopic to the
4-punctured sphere on the other side of the $A$-region, as in
figure \ref{Aregion}). We may thus assume that the surface does
not go through the loops $B_{i,j}$ where $i-j\pmod{n}=2$ (fig.
\ref{complex}(b)). By induction, we may assume that for
$i-j\pmod{n}\leq k$, the torus does not go through $B_{i,j}$,
where $k < n-3$ (fig. \ref{complex}(c-d)).

\begin{figure}[htb]
    \begin{center}
    \epsfxsize=\textwidth
    \epsfbox{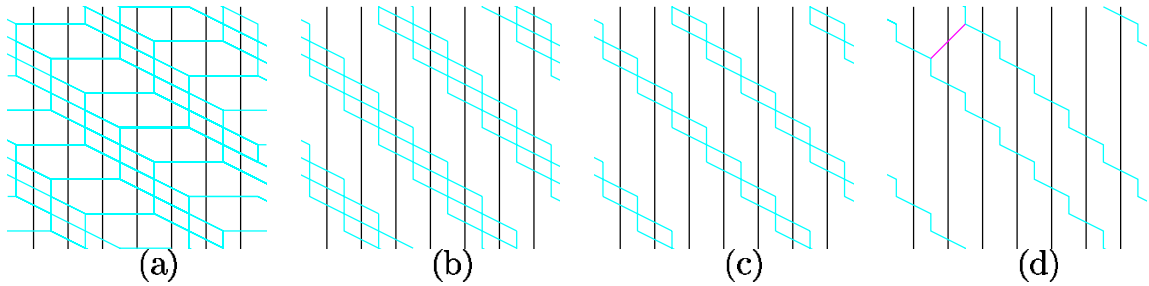}
    \caption{The intersection of the complex with $\g\times I$ for $\S_7$. \label{complex}}
    \end{center}
\end{figure}

\begin{figure}[htb]
    \begin{center}
    \epsfxsize= 3 in
    \epsfbox{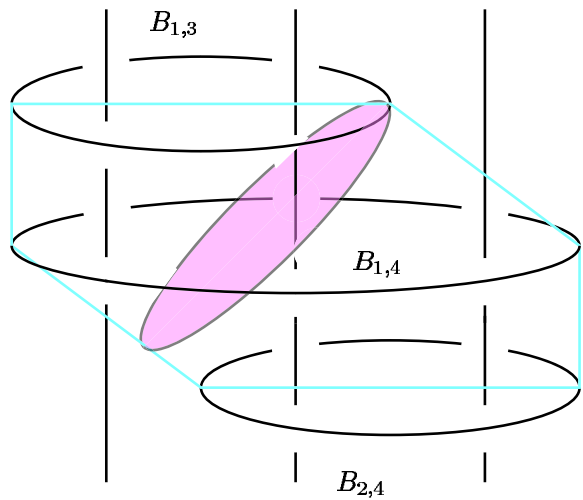}
    \caption{An annular compression of $T$. \label{annuluscompress}}
    \end{center}
\end{figure}

Thus, we may assume that the surface goes through the sequence of
loops $B_{i,j}$, where $i-j\pmod{n}=n-3, n-2$, which gives the
punctured torus $T$. But $T$ is pairwise compressible, which can
be seen by considering the 2-punctured annulus spanning
$\{B_{1,3},B_{1,4},B_{2,4}\}$ (fig. \ref{annuluscompress} gives a
picture of the annular compression). This shows that there is no
pairwise incompressible surface in $T_\psi$ which is isotopic into
$M_C$. If one performs Dehn twists about the curves $B_{i,j}$
avoiding the finite set of Lagrangian subspaces coming from
Hatcher's theorem, one gets a small pure braid. In fact, one may
get a link in $S^3$ by considering $n-1$ strands as lying in a
solid torus, and embedding this solid torus standardly in $S^3$.
\end{proof}

\section{Conclusion}
Given some natural complexity on $3$-manifolds, such as the
minimal number of tetrahedra in a triangulation, or the minimal
number of intersections in a Heegaard diagram, it would be
interesting to understand how the density of small 3-manifolds of
complexity $n$ behaves as $n\longrightarrow \infty$, although this
is probably a very difficult problem. A related question would be
to consider all 2-complexes coming from closed paths in
$\mathcal{P}(\S)^{(1)}$, and to consider the fraction of these for
which the construction given in theorem \ref{smallbraids} produces
manifolds by Dehn twisting along the curves in the pants
decompositions which have no incompressible surface which is
isotopic off of the loops, for a given length of paths in the
pants complex. We conjecture that this fraction goes to 1 as the
length of the paths goes to $\infty$. Indeed one can generalize
the argument of theorem \ref{smallbraids} to show that if a path
goes through a subsequence of pants decompositions as defined in
the proof of the theorem, then any closed incompressible surface
in the complement of the pants curves will compress under Dehn
twisting.

It should be possible to generalize the construction in this paper
to prove the existence of knots in $S^3$ with arbitrarily large
Heegaard genus. It would be interesting to understand this
construction for surfaces with genus. One might be able to prove
the existence of fibred manifolds with fiber of arbitrarily high
genus, such that the fiber is the only closed connected
incompressible surface.


\begin{thebibliography}{10}

\bibitem{Bro01}
{\sc J.~F. Brock}, {\em {Weil-Petersson translation distance and
volumes of
  mapping tori}}, arXiv:math.GT/0109050.

\bibitem{CS99}
{\sc D.~Cooper and M.~Scharlemann}, {\em The structure of a
solvmanifold's
  {H}eegaard splittings}, in Proceedings of 6th G\"okova Geometry-Topology
  Conference, vol.~23, 1999, pp.~1--18.

\bibitem{CJR82}
{\sc M.~Culler, W.~Jaco, and H.~Rubinstein}, {\em Incompressible
surfaces in
  once-punctured torus bundles}, Proc. London Math. Soc. (3), 45 (1982),
  pp.~385--419.

\bibitem{D99}
{\sc N.~Dunfield}, {\em Which small volume hyperbolic 3-manifolds
are haken?}
\newblock talk given at University of Warwick, slides available at
  http://www.math.harvard.edu/~nathand/preprints/slides/haken\_slides.ps, July
  1999.

\bibitem{FH82}
{\sc W.~Floyd and A.~Hatcher}, {\em Incompressible surfaces in
punctured-torus
  bundles}, Topology Appl., 13 (1982), pp.~263--282.

\bibitem{Hak68}
{\sc W.~Haken}, {\em Some results on surfaces in $3$-manifolds},
in Studies in
  Modern Topology, Math. Assoc. Amer. (distributed by Prentice-Hall, Englewood
  Cliffs, N.J.), 1968, pp.~39--98.

\bibitem{HM93}
{\sc J.~Hass and W.~Menasco}, {\em Topologically rigid non-{H}aken
  $3$-manifolds}, J. Austral. Math. Soc. Ser. A, 55 (1993), pp.~60--71.

\bibitem{HLS00}
{\sc A.~Hatcher, P.~Lochak, and L.~Schneps}, {\em On the
{T}eichm\"uller tower
  of mapping class groups}, J. Reine Angew. Math., 521 (2000), pp.~1--24.

\bibitem{HT85}
{\sc A.~Hatcher and W.~Thurston}, {\em Incompressible surfaces in
$2$-bridge
  knot complements}, Invent. Math., 79 (1985), pp.~225--246.

\bibitem{Ha82}
{\sc A.~E. Hatcher}, {\em On the boundary curves of incompressible
surfaces},
  Pacific J. Math., 99 (1982), pp.~373--377.

\bibitem{JS79}
{\sc W.~H. Jaco and P.~B. Shalen}, {\em Seifert fibered spaces in
  $3$-manifolds}, Mem. Amer. Math. Soc., 21 (1979), pp.~viii+192.

\bibitem{Lo92}
{\sc L.~M. Lopez}, {\em Alternating knots and non-{H}aken
$3$-manifolds},
  Topology Appl., 48 (1992), pp.~117--146.

\bibitem{Lo93}
\leavevmode\vrule height 2pt depth -1.6pt width 23pt, {\em Small
knots in
  {S}eifert fibered $3$-manifolds}, Math. Z., 212 (1993), pp.~123--139.

\bibitem{Mar02}
{\sc D.~Margalit}, {\em {The automorphism group of the pants
complex}},
  arXiv:math.GT/0201319.

\bibitem{Mo}
{\sc J.~W. Morgan}, {\em On {T}hurston's uniformization theorem
for
  three-dimensional manifolds}, in The Smith conjecture (New York, 1979),
  Academic Press, Orlando, FL, 1984, pp.~37--125.

\bibitem{MR97}
{\sc Y.~Moriah and H.~Rubinstein}, {\em Heegaard structures of
negatively
  curved $3$-manifolds}, Comm. Anal. Geom., 5 (1997), pp.~375--412.

\bibitem{O84}
{\sc U.~Oertel}, {\em Closed incompressible surfaces in
complements of star
  links}, Pacific J. Math., 111 (1984), pp.~209--230.

\bibitem{Tho}
{\sc A.~Thompson}, {\em Thin position and the recognition problem
for ${S}\sp
  3$}, Math. Res. Lett., 1 (1994), pp.~613--630.

\bibitem{Th}
{\sc W.~P. Thurston}, {\em The geometry and topology of
3-manifolds}.
\newblock Lecture notes from Princeton University, 1978--80.

\bibitem{Th:82}
{\sc W.~P. Thurston}, {\em Three-dimensional manifolds, {K}leinian
groups and
  hyperbolic geometry}, Bull. Amer. Math. Soc. (N.S.), 6 (1982), pp.~357--381.

\bibitem{Wal:78}
{\sc F.~Waldhausen}, {\em Recent results on sufficiently large
$3$-manifolds},
  in Algebraic and geometric topology (Proc. Sympos. Pure Math., Stanford
  Univ., Stanford, Calif., 1976), Part 2, Amer. Math. Soc., Providence, R.I.,
  1978, pp.~21--38.

\end{thebibliography}
\def\cprime{$'$} \def\cprime{$'$}

\end{document}